\documentclass[a4paper,10pt, oneside]{article}
\pdfoutput=1

\usepackage[numbers]{natbib}
\usepackage[utf8]{inputenc}
\usepackage[english]{babel}
\usepackage{amsmath,amsfonts,amssymb, amsthm,bbm}
\usepackage{mathtools}
\usepackage{float}
\usepackage{url}
\usepackage{booktabs}
\usepackage{setspace}
\usepackage{graphicx}
\usepackage[usenames,dvipsnames]{xcolor}
\usepackage{multirow}
\usepackage{stmaryrd}
\usepackage{geometry}
\usepackage{tikz}
\usepackage{pgfplots}
\usepackage[disable]{todonotes}
\usepackage{hyperref}
\usepackage[ruled,vlined,linesnumbered]{algorithm2e}
\SetKwInput{KwPar}{Parameter}

\newtheorem{theorem}{Theorem}%[section]

\newtheorem{proposition}[theorem]{Proposition}

\newtheorem{definition}[theorem]{Definition}

\theoremstyle{definition}

\newtheorem{remark}[theorem]{Remark}

\newcommand{\1}{ \mathbbm{1}}

\newcommand{\cB}{\mathcal{B}}

\newcommand{\cF}{\mathcal{F}}

\newcommand{\cL}{\mathcal{L}}

\newcommand{\cQ}{\mathcal{Q}}

\newcommand{\cT}{\mathcal{T}}

\newcommand{\EE}{\mathbb{E}}

\newcommand{\PP}{\mathbb{P}}
\newcommand{\QQ}{\mathbb{Q}}
\newcommand{\RR}{\mathbb{R}}

\newcommand{\nc}[1]{[#1]}                                   % Crochet normal
\newcommand{\na}[1]{\{#1\}}
\def\va#1{{\boldsymbol{#1}}}

\newcommand{\VA}[1]{\left[\va{#1}\right]}

\def\espe{\EE}

\newcommand{\nesp}[2][]{\espe_{#1}\nc{#2}}           % Esp\'{e}rance normal
\newcommand{\st}{\text{s.t.}}                            % ``s.t.''
\def\eqfinp{\; .}                                           % Point en fin d'\'{e}quation
\def\pfd{\boxtimes}
\DeclareMathOperator{\avar}{AV@R}

\newcommand{\Dfrak}{\mathfrak{D}}

\newcommand{\new}[1]{#1}
\newcommand{\D}{\mathfrak{D}}
\newcommand{\jh}{\hat{\jmath}}

\title{Dual SDDP for risk-averse multistage stochastic programs}
\author{Bernardo Freitas Paulo da Costa,\thanks{EMAp-FGV, Rio de Janeiro, Brazil}
\and Vincent Lecl\`{e}re,\thanks{CERMICS, École des Ponts, Marne-la-Vallée, France}}

\begin{document}
\maketitle

\begin{abstract}
Risk-averse multistage stochastic programs appear in multiple areas and are challenging to solve. 
Stochastic Dual Dynamic Programming (SDDP) is a well-known tool to address such problems under time-independence assumptions.
We show how to derive a dual formulation for these problems and apply an SDDP algorithm, leading to converging and deterministic upper bounds for risk-averse problems.

  \paragraph{Keywords.}
  Stochastic programming, Dynamic programming, SDDP, Risk measures, Duality

  \paragraph{AMS subject classification.}
  90C15, 90C39, 49N15
\end{abstract}

%\tableofcontents

\section{Introduction\label{sec:introduction}}

%!TEX root=../main.tex
%TODO add references.

Multistage stochastic programming is a powerful framework with multiple applications \cite{%wallace2005applications,
gassmann2013stochastic}, \emph{e.g.} in the finance, energy and supply chain sectors.
If the uncertainty is finitely supported, those problems can be seen as large-scale deterministic problems.
When there is more than $4$ or $5$ stages,
the deterministic equivalent is usually too large to be solved directly.
One of the most successful paradigms in this setting
consists in leveraging time-independance assumptions
to derive Bellman equations \cite{bertsekas1995dynamic}.
The Stochastic Dual Dynamic Programming (SDDP) algorithm, and its numerous variants (\cite{pereira1991multi,baucke2017deterministic,zou2019stochastic,ahmed2020stochastic}), consists in using those equations to derive approximations of the \emph{cost-to-go} functions.
It has been successfully used on a number of real-world problems, especially in the field of energy.

While the classical formulation of a multistage program is \emph{risk-neutral}, meaning that we minimize an expected cost,
a large part of the recent litterature sparked by \cite{shapiro2012minimax,philpott2013solving, shapiro2013risk} %TODO : add Dave Morton ?
has been devoted to efficiently introduce \emph{risk aversion} in this framework, in particular inside the SDDP algorithm.
Coherent risk measures \cite{artzner1999coherent} have become a usual tool
to represent risk aversion in stochastic optimization problems.
In multistage stochastic programming,
minimizing a risk measure of the sum of costs leads to time-inconsistency.
% TODO: ref
The easiest way to come up with a time-consistent risk-averse problem is to use composed Markovian risk measures \cite{ruszczynski2010risk}, which, roughly speaking, means replacing the expectation by a risk measure inside the dynamic programming equation.

More precisely, let $(\Omega,\cF, \PP)$ be a probability space,
and $\{\va \omega_t\}_{t\in [T]}$ be a sequence
\new{of finitely supported, independent} random variables
\new{(by convention, boldscript refers to random variables,
normal script to an element of their support,
and equalities between random variables hold almost surely.)}
We consider the following risk-averse multistage linear program (RA-MSLP)
\begin{subequations}
\label{pb:ra-mlsp}
\begin{align}
	\min_{\va x_t, \va y_t} \; & \rho_1\bigg(\va c_1^\top \va y_1 + \rho_{2|\va \omega_1}\Big( \dots + \rho_{T|\va \omega_{[T-1]}}(\va c_T^\top \va y_T)\Big) \bigg) \\
	\st \; & \va A_t \va x_{t} + \va B_t \va x_{t-1} + \va T_t \va y_{t} = \va d_{t} \quad \forall t\in [T] \\
	 &0 \leq \va x_t \new{\leq \bar x_t}, 
   \; 0 \leq \va y_t \new{\leq \bar{y}_t}
   \hspace{0.95cm}
   \forall t \in[T] \label{cst:upperbounds}\\
	 & \va x_t, \va y_t \preceq \va \omega_{[t]}
   \hspace{2.9cm} \forall t\in [T] \label{cst:measurability}
\end{align}
\end{subequations}
where $\va \rho_{t|\omega_{[t]}}$ is a coherent risk measure conditional on the past noises $\va \omega_{[t]}:=\{\va \omega_1, \dots, \va \omega_t\}$, 
all equalities hold almost surely,
and constraint \eqref{cst:measurability} is the \emph{non-anticipativity} constraint,
stating that decisions $\va x_t, \va y_t$
are measurable with respect to $\va \omega_{[t]}:=\na{\va \omega_1, \dots, \va \omega_{t}}$.
\new{Convexity of $\rho_t$ is crucial both for the SDDP algorithm and the duality theory developed here. 
Moreover, in this paper we restrict ourselves to polyhedral risk measures (defined in Section~\ref{sec:duals_poly_rm}) to avoid dealing with technical constraint qualification considerations which would distract the reader.}
\new{Finally, note that, by construction, the nested multistage risk measure used in Problem~\eqref{pb:ra-mlsp} is time-consistent.}

Since $\{\va \omega_t\}_{t\in [T]}$
is a sequence of \emph{independent} random variables,
Dynamic Programming leads to the following recursion:
\begin{align}
  \label{eq:DP_terminal}
  V_{T+1}(x_T) & = 0, \\
  \label{eq:DP_ra_primal}
  V_t(x) & = \!\!
  \begin{array}[t]{rl}
    \min\limits_{\va x_t,\va y_t} \!\!
    & \rho_t\left[ \va c_t^\top \va y_{t} + V_{t+1}(\va x_t)\right] 
    \\
    \textrm{s.t.}
    & \va A_t \va x_t + \va B_t x + \va T_t \va y_t = \va d_t \\
    & 0 \leq \va x_t \new{\leq \bar x_t}, 
    \; 0 \leq \va y_t \new{\leq \bar{y}_t}
  \end{array}
\end{align}
where the value of Problem~\eqref{pb:ra-mlsp} is given by $V_1(x_0)$.

The classical SDDP algorithm builds outer approximations
of the cost-to-go functions $V_t$,
leading to exact lower bounds on the problem.
In a risk-neutral framework, upper bounds can be estimated via Monte Carlo sampling.
Unfortunately, it is unclear how to extend such statistical methods
to the risk-averse setting~\cite{shapiro2013risk}.
Instead of statistical upper bounds, one can use exact upper bounds:
Through backward recursion (\cite{philpott2013solving});
by maintaining upper and lower bounds for all value functions (\cite{baucke2017deterministic,downward2020stochastic});
or using Fenchel duality (\cite{leclere2020exact,guigues2019duality}).
\new{Up to now, the first approach has not been used
to compute improving upper bounds along SDDP iterations,
while the second approach relies on a \emph{problem-child} node selection method.
Finally, the last approach was developed only
in a risk-neutral setting.
The aim of this work is to adapt the latter approach to a risk-averse setting.
By dualizing the extensive formulation of the risk-averse MLSP problem, and recognizing a time-decomposition, we obtain a Bellman recursion on which SDDP can be applied, yielding converging exact upper bounds.}

\medskip

\paragraph{Contributions}
In this paper we i) derive a dual formulation of RA-MLSP with polyhedral risk measure;
ii) show that it is time-decomposable and solvable through SDDP, yielding exact upper bounds of the original problem;
iii) link the value function of the dual formulation with the co-perspective of the primal value function;
and
iv) illustrate the approach with numerical results.

% TODO: notation [J] = {1, 2, ..., J} and \omega_{[t]} = (omega_1, ..., \omega_t)

% vim:set spelllang=en:

% \section{An illustrative example}
% \label{sec:example}
% \input{subfiles/three_stage_example}

\section{Time decomposition of the dual of a risk averse MSLP}
\label{sec:decomposition}

\subsection{Risk-averse duals with AV@R}
\label{sec:examples_of_duals}
%!TEX root=../main.tex

%\subsection{Dual formulation for the mean-AV@R risk measure}

We start by showing how to build the dual problem in a very specific setting:
for a single step of the recursion,
with no upper bounds on $\va{x}_t$ and $\va{y}_t$,
and when the risk measure $\rho$ is
a convex combination of the mean and the $\alpha$-AV@R,
given by, for \new{$\alpha \in (0,1)$ and $\beta \in [0,1]$},
\begin{equation}
  \label{mean_avar}
  \rho\VA{\theta} := \beta \nesp{\va{\theta}} + (1 - \beta) \avar_\alpha[\va{\theta}] \eqfinp
\end{equation}
This risk measure assumes an underlying probability for the scenarios,
with respect to which one calculates the expectation and the AV@R.
% Moreover, it is decreasing in $\alpha$ and $\beta$.
The risk measures we employ in the example in section~\ref{sec:numerical}
will be of this class.

We rewrite equation~\eqref{eq:DP_ra_primal}
using the Rockafellar-Uryasev representation of $\avar$,
with $\va{\theta}$ as epigraphical variables for the scenario costs.
\new{For simplicity, we represent a random variable as a vector in $\RR^J$,
denoted with bold letters such as $\va{x} = (x_1, \dots x_J)$,
and the expectation $\EE[\va{x}]$ is the sum $\sum p_j x_j$.}
So, the value of $V_t(x_{t-1})$ is given by:
\begin{equation}
  \label{DP_mean_avar}
  \begin{array}[t]{rll}
    \inf\limits_{\va{x}, \va{y}; q, \va{\theta}, \va{u}}
      & \multicolumn{2}{l}{\beta \nesp{\va{\theta}} + (1 - \beta) \left[q + \frac{1}{\alpha} \nesp{\va{u}} \right]} \\[.5ex]
    \textrm{s.t.}
      & q + {u}_j \geq {\theta}_j & \forall j\in[J] \\[.3ex]
      & {\theta}_j \geq c_j^\top y_j + V_t(x_j) & \forall j\in[J] \\[.3ex]
      & A_j x_j + B_j x_{t-1} + T_j y_j = d_j & \forall j\in[J] \\[.3ex]
      & x_j, y_j, u_j \geq 0 &  \forall j\in[J]
  \end{array}
\end{equation}

We define dual multipliers for every constraint:
in order, $\va{\delta}$, $\va{\gamma}$, $\va{\lambda}$,
$\va{\mu}$, $\va{\nu}$, and~$\va{\eta}$.
With the expectation inner product, this yields the Lagrangian:
\begin{multline*}
  (1 - \beta) q + \beta \nesp{\va{\theta}} + (1-\beta)/\alpha \cdot \nesp{\va{u}}
  + \nesp{ \va{\gamma} ( \va{c}^\top \va{y} + V(\va{x}) - \va{\theta} )
         + \va{\delta} ( \va{\theta} - q - \va{u} ) } \\
  \quad + \nesp{ \va{\lambda}^\top ( \va{A} \va{x} + \va{B} x_{t-1} + \va{T} \va{y} - \va{d} ) }
  - \nesp{ \va{\mu}^\top \va{x} + \va{\nu}^\top \va{y} + \va{\eta} \va{u} }
\end{multline*}
Eliminating the multipliers $\va{\nu}$ and~$\va{\eta}$,
we obtain the dual problem
\begin{equation}
  \label{dual_mean_avar}
  \begin{array}[t]{rll}
    \sup\limits_{\va{\lambda}, \va{\gamma}, \va{\delta}, \va{\mu}}
    & \multicolumn{2}{l}{\begin{array}[t]{l}
      % TODO: verify sign here: + or - inf?
      \mathbb{E} \bigg[ \va{\lambda}^\top(\va{B} x_{t-1} - \va{d}) +
        \inf\limits_{\va{x}} \left[ (\va{A}^\top \va{\lambda} - \va{\mu})^\top \va{x} + \va{\gamma}V_t(\va{x}) \right] \bigg]
    \end{array}} \\[6ex]
    \textrm{s.t.}
    & \nesp{\va{\delta}} = (1 - \beta) \\
    & 0 \leq {\delta}_j \leq \frac{1 - \beta}{\alpha} & \forall j\in[J] \\
    & {\gamma}_j = \beta + {\delta}_j & \forall j\in[J] \\
    & \gamma_j {c}_j + T_j^\top \lambda_j \geq 0 \quad & \forall j\in[J] \\
    & {\mu}_j \geq 0 & \forall j\in[J]
  \end{array}
\end{equation}

Observe that the variable $\va{\gamma}$
represents the ``change-of-measure'' implied by the
mean-$\avar$ combination~\cite{shapiro2009lectures}.
% TODO? the change of measure can be described in a constructive way (cf., [15, Remarks 24-25, pp.314-315]).
Indeed, $\va{\gamma}$ is at least $\beta \leq 1$,
and some events will have an increased contribution,
up to $\frac{1-\beta}{\alpha}$, so that $\nesp{\va{\gamma}} = 1$.

% vim:set spelllang=en:

\subsection{Polyhedral risk measures and duality}
\label{sec:duals_poly_rm}
To extend the previous approach to more general risk measures,
we adopt a distributionally robust point of view.
We consider a \emph{polyhedral risk measure $\rho$}, that is,
a coherent risk measure of the form
\begin{equation}
  \label{poly_rm}
  \rho : \va t \mapsto \sup_{\QQ \in \cQ} \nesp[\QQ]{\va{t}}
  = \max_{k\in[K]} \{ \nesp[\QQ^k]{\va{t}} \},
\end{equation} 
where 
$\cQ = \mathrm{conv}(\{\QQ^k\}_{k\in[K]})$.
Polyhedral risk measures can be either chosen as interpretable risk-measures
(e.g. $\avar$ in a finite setting)
or as the worst case among a set of probabilities estimated by various experts.
Since we don't assume a reference probability,
we resort to describing the extremal risk measures,
which may be very numerous.
This also changes the interpretation of the dual variables $\va{\gamma}$:
now they correspond to supporting probabilities,
instead of a change-of-measure.

% TODO: assume finitely many scenarios from the outset?
We denote the elements of the support of $\va \omega$ by
$\omega_1, \dots, \omega_J$,
and let $q^k_j := \QQ^k[\va \omega = \omega_j]$.
Now, $V_t(x_{t-1})$ is given by:
\begin{alignat}{4}
  \inf\limits_{\va{x}, \va{y}; z, \va{\theta}} \quad
    % Or: \inf\limits_{\va{x}_1, \va{y}_1 \geq 0; z_0, \va{t}_1} ??
      & z \label{pb:ra_lbo_poly} \\[.5ex]
    \textrm{s.t.} \quad
      & z \geq \sum_{j \in [J]} q^k_j \theta_j & \forall k  &\;& [\phi_k]\nonumber\\[.5ex]
      & \theta_j \geq c_j^\top y_j + V_{t+1}(x_j) & \forall j &\;& [\gamma_j]\nonumber\\[.3ex]
      & A_j x_j + B_j x_{t-1} + T_j y_j = d_j \quad & \forall j  &\;& [\lambda_j]\nonumber\\[.3ex]
      & 0 \leq x_j \leq \bar{x}_t & \forall j &\;& [\mu_j, \zeta_j] \nonumber \\
      & 0 \leq y_j \leq \bar{y}_t & \forall j &\;& [\nu_j, \xi_j] \nonumber
\end{alignat}

% \begin{equation}
%   \label{pb:ra_lbo_poly}
%   \begin{array}[t]{rllr}
%     \inf\limits_{\va{x}, \va{y}; z, \va{\theta}}
%     % Or: \inf\limits_{\va{x}_1, \va{y}_1 \geq 0; z_0, \va{t}_1} ??
%       & z \\[.5ex]
%     \textrm{s.t.}
%       & z \geq \sum_{j \in [J]} q_{k,j} \theta_j & \forall k  & [\phi_k]\\[.5ex]
%       & \theta_j \geq c_j^\top y_j + V(x_j) & \forall j & [\gamma_j]\\[.3ex]
%       & A_j x_j + B_j x_{t-1} + T_j y_j = d_j & \forall j  & [\lambda_j]\\[.3ex]
%       & x_j, y_j \geq 0 &  \forall j & [\mu_j, \nu_j]
%   \end{array}
% \end{equation}

Proceeding analogously to the AV@R case above,
we introduce dual multipliers
%$\phi_k$, $\gamma_j$, $\lambda_j$, $\mu_j$, $\nu_j$
as indicated in the brackets,
and obtain the following dual problem
\begin{alignat}{5}
  \sup\limits_{\substack{\phi_k, \gamma_j, \lambda_j, \\ \mu_j,\new{\zeta_j,\xi_j}}} &
      \sum\limits_{j\in[J]} \Big[ \lambda_j^\top \left(B_j x_{t-1} - d_j\right) \new{- \bar{x}_t \zeta_j } \new{- \bar{y}_t \xi_j} \label{lpdual_ra2} \\
      & \qquad {} + \inf\limits_{x_j} (A_j^\top\lambda_j - \mu_j \new{+ \zeta_j})^\top x_j + \gamma_j V_{t+1}(x_j) \Big] \nonumber\\
    \textrm{s.t.}\quad
    & \sum_k \phi_k = 1 , \qquad \phi_k \geq 0, \nonumber\\
    & \sum_k \phi_k q^k_j = \gamma_j \geq 0 & \forall j\nonumber\\
    & \gamma_j c_j + T_j^\top \lambda_j  \new{+\xi_j}\geq 0 & \forall j \nonumber\\
    & \mu_j, \zeta_j, \xi_j \geq 0 & \forall j. \nonumber
\end{alignat}

% \begin{equation}
%   \label{lpdual_ra2}
%   \begin{array}[t]{rll}
%     \sup\limits_{\phi_k, \gamma_j, \lambda_j, \mu_j} &
%     \multicolumn{2}{l}{\begin{array}[t]{l}
%       \sum_{j\in[J]} \Big[ \lambda_j^\top \left(B_j x_0 - d_j\right) + {} \\
%       \ \ \inf\limits_{x_j} \gamma_j V(x_j) + (A_j^\top\lambda_j - \mu_j)^\top x_j \Big]
%     \end{array}} \\[6ex]
%     \textrm{s.t.}
%     & \sum_k \phi_k = 1 , \qquad \phi_k \geq 0, \\
%     & \sum_k \phi_k q_{k,j} = \gamma_j \geq 0 & \forall j\\
%     & \gamma_j c_j + T_j^\top \lambda_j \geq 0 & \forall j \\
%     & \mu_j \geq 0 & \forall j.
%   \end{array}
% \end{equation}
The constraints on $\phi_k$ are equivalent to describing the vector of $\gamma_j$'s
as a convex combination of the extreme probabilities $\QQ^k$.
Therefore, one can rewrite problem~\eqref{lpdual_ra2} to include the constraint
$\{\gamma_j\}_{j \in [J]} \in \cQ$
instead of the first two lines.
This shows that the variables $\gamma_j$ correspond to one
supporting probability of the risk measure $\rho$.
In particular, if a given scenario is \emph{effective},
in the sense of \cite{rahimian2019identifying},
then there exists an optimal $\va\gamma$ which charges this scenario.

Moreover, the last two constraints here
correspond exactly to the last two in problem~\eqref{dual_mean_avar},
which emphasizes the similarity
between~\eqref{lpdual_ra2} and~\eqref{dual_mean_avar}.

% Moreover,
% \begin{align*}
%   & \inf_{x_j} \gamma V(x_j) + (A_j^\top \lambda_j - \mu_j)x_j \\
%   & \quad = -  \sup_{x_j} - \gamma V(x_j) - (A_j^\top \lambda_j - \mu_j)^\top x_j \\
%   & \quad = -  V^\pfd(\mu_j - A_j^\top \lambda_j, \gamma)
% \end{align*}
% Therefore, the risk-averse linear Bellman operator $\cB$ is also given by
% \begin{equation}
%   \label{risky_lbo}
%   \cB(V)(x_0)
%   =
%   \begin{array}[t]{rll}
%     \sup\limits_{ {\gamma}, {\lambda}, {\mu}}
%     & \multicolumn{2}{l}{\begin{array}[t]{l}
%       \sum_{j\in[J]} \Big[ \lambda_j^\top \left(B_j x_0 - d_j\right) - {} \\
%         V^\pfd(\mu_j - A_j^\top \lambda_j, \gamma_j) \Big]
%     \end{array}} \\[6ex]
%     \textrm{s.t.}
%     & \gamma \in \cQ  \\
%     & \gamma_j c_j + T_j^\top \lambda_j \geq 0 & \forall j  \\
%     &  \mu_j \geq 0 .
%   \end{array}
% \end{equation}

% vim: set spelllang=en:

%!TEX root=../main.tex
\subsection{Multistage risk averse problem duality}

We now extend the duality to the full multistage problem.
In the stagewise independent setting, we let
$\Omega_t$ be the set of all possible realizations of $\omega_t$,
and the risk measure $\rho_t$ is defined by
$\rho_t = \sup_{\QQ \in \cQ_t} \EE_\QQ [ \cdot]$,
for a polyhedral subset $\cQ_t$ of probability measures on $\Omega_t$.
The tree $\cT$ describing the stochastic process is such that
each node $n$ of depth $t$ is associated with a possible value of
$\va\omega_{[t]} = (\va\omega_1, \ldots \va\omega_t)$.
For any node $n$, the set of its children is denoted by $C_n$,
and $\cL$ is the set of leaves of $\cT$.

In the spirit of the previous section,
we introduce variables $z_n$ to stand for the risk-adjusted value
of our problem starting from node $n$,
and $\theta_m$ represents the cost-to-go following the branch of node $m \in C_n$.
To reduce notational burden, we assume that, for all $t$,
$\rho_t = \rho$.
Then, the risk averse problem~\eqref{pb:ra-mlsp},
with value $V_{n_0}(\tilde x_{n_0})$,
can be written as the following linear program:
\begin{alignat}{4}
	\min & \quad z_0  \label{pb:RAMSLP-extended}\\
  \text{s.t.}
		& \sum\limits_{m \in C_n} q_m^k \theta_m  \leq z_n
		  & \forall n, \forall k \in [K]
		  &\ & [\Phi_n^k] \nonumber\\
		& c_m^\top y_m + z_m  \leq  \theta_m
		  & \ \ \forall m \in \cT \backslash \{n_0\}
		  &\ & [\gamma_m] \nonumber\\
		& A_m x_m + B_m \tilde{x}_n + T_m y_m = d_m \quad
		  & \forall n, \forall m \in C_n
		  &\ & [\lambda_m] \nonumber\\
		& z_\ell = 0
		  & \forall \ell \in \cL
		  &\ & [\eta_\ell]\nonumber\\
		& x_m = \tilde{x}_m
		  & \forall m \in \cT \backslash \{n_0\}
		  &\ & [\pi_m]\nonumber\\
		& 0 \leq \tilde x_m \leq \bar{x}_m
		  & \forall m \in \cT \backslash \{n_0\}
		  &\ & [\mu_m, \zeta_m]\nonumber\\
		& 0 \leq y_m \leq \bar{y}_m
      & \forall m \in \cT \backslash \{n_0\}
      &\ & [\mu_m,\xi_m]
		  \nonumber
\end{alignat}
	% \begin{equation}
	% \label{pb:RAMSLP-extended}
	%   \def\arraystretch{1.25}
	%   \begin{array}{rlll}
	% 	\min  & z_0 \\
	% 	& \sum\limits_{m \in C_n} q_m^k \theta_m  \leq z_n
	% 	  & \forall n, \forall k \in [K]
	% 	  & [\Phi_n^k] \\
	% 	& c_m^\top y_m + z_m  \leq  \theta_m
	% 	  & \forall m \in \cT \backslash \{n_0\}
	% 	  & [\gamma_m] \\
	% 	& A_m x_m + B_m \tilde{x}_n \\
	% 		& \quad
	% 	+ T_m y_m = d_m
	% 	  & \forall n, \forall m \in C_n
	% 	  & [\lambda_m] \\
	% 	& x_n \geq 0, y_n \geq 0
	% 	  & \forall n \in \cT \backslash \{n_0\}
	% 	  & [\mu_n, \nu_n]\\
	% 	& z_\ell = 0
	% 	  & \forall \ell \in \cL
	% 	  & [\eta_\ell]\\
	% 	& x_n = \tilde{x}_n
	% 	  & \forall n 
	% 	  & [\pi_n]\\
	% 	  &x_n \leq \bar{x}_n,\quad y_n \leq \bar{y}_n &\forall n & [\zeta_n,\xi_n]
	%   \end{array}
	% \end{equation}
where, when unspecified, $\forall n$ stands for $\forall n \in \cT \backslash \cL$,
$\tilde x_{n_0}$ is a parameter and not a variable,
and we add the equalities $x_m = \tilde{x}_m$ to highlight the time dynamics.

Defining $\gamma_{n_0} = 1$,
the linear programming dual of problem~\eqref{pb:RAMSLP-extended} is
\begin{alignat*}{3}
  \sup_{\Phi,\gamma, \pi, \lambda} \ 
  % & \pi_{n_0}^\top \tilde{x}_{n_0}- \sum_{m \in \cT \backslash \{n_0\} } \lambda_m^\top d_m \\
  & \multispan3{$\displaystyle \pi_{n_0}^\top \tilde{x}_{n_0} -
                \sum_{m} \lambda_m^\top d_m + \bar{x}_m^\top \zeta_m + \bar{y}_m^\top \xi_m$\hfil} \\
  \text{s.t.} \ 
  & \sum_{k\in [K]} \Phi_n^k = \gamma_n
    &\ & \forall n & [z_n] \\
	& \sum_{k\in [K]} \Phi^k_n q^k_m = \gamma_m \geq 0
    &\ & \forall n, \forall m \in C_n \ & [\theta_m] \\
	& \pi_{n_0} = \sum_{m \in C_{n_0}}  B_m^\top \lambda_m \\
	& \pi_n \leq \zeta_n + \sum_{m\in C_n} B_m^\top \lambda_m
    &\ & \forall n \in \cT \backslash \{n_0\} \ & [\tilde{x}_m] \\
	& \pi_m + A_m^\top \lambda_m = 0
    &\ & \forall m & [x_m] \\
	& \gamma_m c_m + T_m^\top \lambda_m  + \xi_m \geq 0
    &\ & \forall m & [y_m] \\
	& \Phi_n^k \geq 0
    &\ & \forall n, \forall k \in [K]\\
	& \zeta_m \geq 0, \xi_m \geq 0
    &\ & \forall m
\end{alignat*}
where we keep $\forall n$ to imply $n \in \cT \backslash \cL$ as above,
and unspecified $\forall m$, $\sum_m$ range over $m \in \cT\backslash \{n_0\}$.
%\begin{strip}
% \begin{align*}
% 	D_{n_0}(\tilde x_{n_0})=\sup_{\Phi,\gamma , \lambda}  \qquad& \pi_{n_0}^\top \tilde{x}_{n_0}- \sum_{m \in \cT \backslash \{n_0\} } \lambda_m^\top d_m \\
% 	%& \sum_{k \in [K]} \Phi^k_0 = 1 & [z_0] \\
% 	& \sum_{k\in K} \Phi_n^k = \gamma_n & \forall n \in \cT \backslash \cL \qquad[z_n] \\
% 	& \sum_{k\in [K]} \Phi^k_n q^k_m = \gamma_m & \forall n \in \cT\backslash \cL, \forall m \in C_n  \qquad [\theta_m]\\
% 	& \gamma_m c_m + T_m^\top \lambda_m \geq 0 & \forall m \in \cT\backslash \{n_0\} \qquad[y_m] \\
% 	& 
% 	\sum_{m\in C_n} B_m^\top \lambda_m = \pi_n & \forall n \in \cT \backslash {\cL} \qquad [\tilde{x}_n] \\
% 	& A_m^\top \lambda_m  + \pi_m \geq 0 &
% 	\forall m \in \cT \backslash \{n_0\} \qquad [x_m]
% 	\\
% 	& \Phi_n^k \geq 0, \gamma_n \geq 0 & \forall n \in \cT \ .
% \end{align*}
%\end{strip}

Note that $\Phi_n^k$ can be seen as barycentric coordinates of the extreme points of $\cQ$.
Thus, the first two constraints can be more compactly written as $(\gamma_m)_{m \in C_n} \in \gamma_n\cQ$.

By backward recursion,
this problem can be solved through the following recursive equations, where,
for all leaves $\ell \in \cL$,
$D_\ell(\pi_\ell, \gamma_\ell) = - \bar{x}_\ell^\top \max\{\pi_\ell, 0\}$,
and for all nodes $n \in \cT\backslash \cL$,
$D_n(\pi_n,\gamma_n)$ is given as the value of
\begin{align}
  \sup_{\substack{\pi_m, \gamma_m, \lambda_m\\ \zeta_n, \xi_m \geq 0}} \ &
   \1_{\{n = n_0\}} \pi_{n_0}^\top \tilde x_{n_0} - \bar{x}_n^\top \zeta_n + \label{eq:Dn}\\
  & \ \sum_{m \in C_n} - \lambda_m^\top d_m - \bar{y}_m^\top \xi_m+ D_m(\pi_m,\gamma_m)  \nonumber\\
  \text{s.t.} \quad	& (\gamma_m)_{m \in C_n} \in \gamma_n \cQ \nonumber\\
  & \zeta_n + \sum_{m \in C_n} B_m^\top \lambda_m \geq \pi_n \nonumber\\
  & \pi_m + A_m^\top \lambda_m = 0, \qquad\qquad\,\ \forall m \in C_n \nonumber\\
  & \gamma_m c_m + T_m^\top \lambda_m  + \xi_m \geq 0, \quad \forall m \in C_n \nonumber
\end{align}

By the independence assumption, a backward induction shows that $D_n = D_{n'}$ for all nodes $n$ and $n'$ of the same depth.
% TODO: check T or T+1
Thus, defining $D_T(\pi_T, \gamma_T) = - \bar{x}_T^\top \max\{\pi_T, 0\}$,
we obtain the following recursion for the dual value functions:
\begin{align}
  & D_t(\pi_{t},\gamma_{t}) = \label{eq:dualRCR} \\
  & \begin{aligned}
    \sup_{\substack{\zeta, \gamma_j, \\ \lambda_j, \pi_j, \xi_j}}
    & - \bar{x}_{t}^\top \zeta
      + \sum_{j \in [J_t]} \Big[ - d_j^\top \lambda_j  - \bar{y}_{t+1}^\top \xi_j
                                 + D_{t+1}(\pi_j,\gamma_j) \Big] \\
  \text{s.t.} \quad
	  & (\gamma_j)_{j \in [J_t]} \in \gamma_t \cQ  \\
	  & \zeta + \sum_{j \in [J_t]} B_j^\top \lambda_j \geq \pi_t \\
	  & \pi_j + A_j^\top \lambda_j = 0 & \forall j \in [J_t] \\
	  & \gamma_j c_j + T_j^\top \lambda_j + \xi_j\geq 0 & \forall j \in [J_t] \\
	  & \xi_j \geq 0, \quad \zeta \geq 0
  \end{aligned}\nonumber
\end{align}
This decomposition satisfies the RCR conditions.
Indeed, for every $\pi_t$ and every $\gamma_t \geq 0$,
any $\gamma \in \gamma_t \cQ$ and $\lambda = 0$ are admissible, using slack variable $\zeta$ as needed.
Then, $\pi_j$ \new{are given by the $\pi_j + A_j^\top \lambda_j = 0$},
and the remaining constraints can be adjusted using $\xi_j$.

%   \new{The explicit upper bound constraints on the state $\va x_t$ and control $\va y_t$ are often reasonable from a modelling point of view, and usefull to ensure relatively complete recourse in the dual (which is not otherwise guaranteed, see \cite{guigues2019duality}).}

% One of the difficulties of the dual formulation is that even if relatively complete recourse (RCR) is ensured in the primal, it is not guaranteed in the dual (see for example \cite{guigues2019duality}).
% To deal with this problem we can either incorporate feasibility cuts in the algorithm, 
% or add a penalization scheme to ensure RCR.
% We follow here a slightly different path, equivalent to the penalization approach:
% we assume that we can bound the state and control variables a priori by constants $\bar{x}_n,\bar{y}_n$, 
% and explicitly add those constraints, yielding new multipliers.

\begin{remark}
	Relatively complete recourse in a dual formulation is not guaranteed (see for example \cite{guigues2019duality}).
	In our setting, the explicit upper bounds of~\eqref{cst:upperbounds} ensure RCR.
	The existence of such upper bounds is equivalent to the existence of exact penalization coefficients in the dual, which
	is the tool used in \cite{guigues2019duality} to deal with this difficulty. 
	Alternatively, we could incorporate feasibility cuts in the algorithm.
\end{remark}

\subsection{Bounding the dual state}
\label{sub:bounding_state}

With our boundedness assumption,
we have relatively complete recourse \new{in the dual}.
To prove convergence, we still need to ensure that
the \new{dual} state remains bounded.

By assumption, we know that there exists an optimal primal solution.
Further, by linear programming duality, we know that there exists an optimal dual solution.
The marginal interpretation of the Lagrange multiplier $\pi$
(see Problem~\eqref{pb:RAMSLP-extended}) states that,
for each node, the optimal dual $\pi_n$ is a subgradient of the primal value function \emph{for $\gamma_n = 1$}. 
In particular, $\pi_n/\gamma_n$ can be bounded by the Lipschitz constant of the primal value function $V_n$.
\new{In the independent setting,}
assuming that $V_t$ is $L_t$-Lipschitz continuous on its domain,
we can add the constraint $|\pi_j| \leq \gamma_j L_{t+1}$ to~\eqref{eq:dualRCR}
for each $j$, without changing its value.
This method is similar to the compactification process through Lipschitz-regularization used in~\cite{leclere2020exact}.

% TODO: comment on the difference with Shapiro-Guigues?

Therefore, we use the compactified \new{recursion}
presented in \eqref{projective_dual_lbo_compacted}.
\new{Since it has RCR and bounded states,}
the SDDP algorithm on this recursion converges.
This is illustrated in section~\ref{sec:examples}.

\begin{equation}
    \label{projective_dual_lbo_compacted}
    D_t(\pi_t, \gamma_t)
    =  \begin{array}[t]{rll} 
  \sup\limits_{\zeta, \gamma_j, \lambda_j, \pi_j, \xi_j} &
  - \bar{x}_t^\top \zeta +
    \sum\limits_{j\in[J]} - d_j^\top \lambda_j
      - \bar{y}_{t+1}^\top \xi_j + D_{t+1}(\pi_j, \gamma_j) \\[.3ex]
    \textrm{s.t.} \quad
    & \gamma \in \gamma_t \cQ  \\
    & \zeta + \sum_j B_j^\top\lambda_j \geq \pi_t \\
    & \pi_j + A_j^\top \lambda_j = 0  & \forall j \in [J_t] \\
    & \gamma_j c_j + T_j^\top \lambda_j + \xi_j \geq 0 & \forall j \in [J_t] \\
    & |\pi_j|\leq \gamma_j L_{t+1} & \forall j \in [J_t] \\
    & \zeta \geq 0, \xi_j \geq 0 \\

  \end{array}
\end{equation}

% vim:set spelllang=en:

\section{Dual risk averse Bellman operator}
\label{sec:dual_bellman_operator}

We introduce convex analysis tools that shed new light on the link between the primal and dual value functions given in Section~\ref{sec:decomposition}.
% More precisely we show that the dual value function $D_t$ introduced in~\eqref{eq:D_t} is the negative of coperspective of the primal value function $V_t$.

\subsection{Homogeneous Fenchel duality}
%!TEX root=../main.tex

Let $f: \RR^n \to (-\infty, \infty]$ be a proper lower semicontinuous convex function. 
Recall (see \cite{combettes2018perspective} for more details) that
the \emph{perspective function} of $f$, denoted $\tilde{f}$,
is a convex, lower-semicontinuous function of $\RR^{n+1}$,
such that $\tilde{f}(x,\gamma) = \gamma f(x/\gamma)$ for any positive number $\gamma$.
% \begin{equation}
%   \label{eq:def_perspective}
%   \begin{aligned}[t]
%     \tilde f : \RR \times \RR^n & \to (-\infty,+\infty] \\
%                      (\gamma,x) & \mapsto
%   \begin{cases}
%     \gamma f(x/\gamma), & \text{ if } \gamma >0 \\
%     \rec(f)(x), & \text{ if } \gamma =0 \\
%     +\infty & \text{otherwise}
%   \end{cases}
%   \end{aligned}
% \end{equation}
% where $\rec(f)$ is the \emph{recession function} of $f$ defined as
% \begin{equation}
%   \label{eq:def_recession}
%   \rec(f) : \RR^n \to ]-\infty,+\infty]:
%   x \mapsto \lim_{t \to +\infty} \frac{f(z+tx)}{t} \; ,
% \end{equation}
% where $z$ can be chosen as any point in the domain of $f$.
% Both the recession and the perspective function of proper lower semicontinuous convex function are proper lower semicontinuous convex function as well.

Recall that the Fenchel conjugate of $f$ is 
\begin{equation}
\label{eq:fenchel}
f^\star : \RR^n \to \overline \RR : \psi \mapsto
\sup_{x \in \RR^n} \psi^\top x - f(x).
\end{equation}

\medskip

Inspired by the recurrences in~\eqref{dual_mean_avar} and~\eqref{lpdual_ra2},
we introduce the \emph{coperspective function}:
\begin{definition}
Let $f: \RR^n \to \overline \RR$.
The coperspective of $f$ is the perspective of the Fenchel conjugate, that is $(f^\star)^\sim$, that we denote $f^\pfd$.
In particular, for $\psi \in \RR^n$ and $\gamma \in \RR_{++}$,
we have
\begin{equation}
  \label{pdf}
  f^\pfd(\psi, \gamma) := \sup_{x \in \RR^n} \psi^\top x - \gamma f(x).
\end{equation}
\end{definition}
% We (probably) need to extend the definition for $\tau = 0$
% so that the function is lower semi continuous;
% the current form is +\infty for all $\pi \neq 0$...

\begin{remark}
The coperspective is jointly convex in $(\psi, \gamma)$,
lower semicontinuous, and a positively homogeneous function of degree $1$:
for all $t>0$,
\[
  f^\pfd(t \cdot \psi, t \cdot \gamma)
    = t \cdot f^\pfd(\psi, \gamma) .
\]
\end{remark}

\begin{remark}
\label{rk:beta0}
Cuts for a convex function and its perspective are essentially equivalent.
If $f(x) \geq f(x_0) + g^\top(x - x_0) = \theta + g^\top x$,
then
\begin{align*}
  \tilde{f}(x,\gamma) = \gamma \cdot f(x/\gamma)
         & \geq \gamma f(x_0) + \gamma g^\top( x/\gamma - x_0) \\
         & \geq \gamma f(x_0) + g^\top( x - \gamma \cdot x_0) \\
         & \geq \theta \cdot \gamma + g^\top x
\end{align*}
Similarly, if $\tilde f(x,\gamma) \geq \theta \cdot \gamma + g^\top x + \beta$,
then $f(x) \geq g^\top x + \theta + \beta$.
Note that if the cut for $\tilde f$ is exact we can assume $\beta = 0$.
\end{remark}

% vim:set spelllang=en:

\subsection{Duality and conjugate value functions}

Consider a polyhedral risk measure $\rho$ and the associated 
risk-averse Bellman operator $\cB$  that, to any cost-to-go function $V$
and initial state $ x_{t-1}$ associates the value of Problem~\eqref{pb:ra_lbo_poly}.

The coperspective of $\cB(V)$ can be calculated using~\eqref{lpdual_ra2}.
Leveraging positive homogeneity,
for $\psi_0 \in  \RR^n$ and $\gamma_0 > 0$,
we get that $\cB(V)^\pfd(\psi_0, \gamma_0)$ is given by
\begin{align}
   \label{fenchel_dual_ra_hom}
    \sup_{x_0} \ \psi_0^\top x_0 \\
    + 
  \inf\limits_{\substack{{\gamma}, {\lambda}, {\mu}\\\new{\zeta, \xi}}} \ &
    \sum_{j\in[J]} \lambda_j^\top \left(d_j - B_j x_0 \right) + \new{\xi_j^\top \bar{y}_{t+1}}
      + \zeta_j^\top \bar{x}_{t+1} + V^\pfd(\mu_j - A_j^\top \lambda_j - \zeta_j, \gamma_j) \nonumber \\[.3ex]
    \textrm{s.t.} \ 
    & \gamma \in \gamma_0 \cQ \nonumber \\
    & \gamma_j c_j \new{+ \xi_j} + T_j^\top \lambda_j \geq 0 \qquad\quad \forall j \nonumber \\
    &  \mu_j \new{, \zeta_j, \xi_j} \geq 0 \qquad\qquad\qquad\ \ \, \forall j. \nonumber
\end{align}

Note that, if $V$ is polyhedral, so are its Fenchel dual and its perspective.
Thus, by linear programming duality,
we can interchange $\sup$ and $\inf$ to obtain
%\begin{strip}
\begin{align}
   \label{fenchel_dual_ra2}
  [\cB(V)]^\pfd(\psi_0, \gamma_0)
  & = \begin{array}[t]{rll}
  \inf\limits_{\substack{{\gamma}, {\lambda}\\\new{\mu, \zeta, \xi}}}
    & \sum\limits_{j\in[J]} \lambda_j^\top d_j + \new{\xi_j^\top \bar{y}_{t+1} + \zeta_j^\top \bar{x}_{t+1} +} V^\pfd(\psi_j, \gamma_j)
        \\[.3ex]
    \textrm{s.t.}
    & \sum_j B_j^\top \lambda_j = \psi_0 \\
    & \gamma \in \gamma_0 \cQ  \\
    & \gamma_j c_j + \xi_j + T_j^\top \lambda_j \geq 0 \quad \forall j  \\
    & \psi_j = \mu_j - A_j^\top \lambda_j - \zeta_j    \quad \forall j .
  \end{array}
\end{align}
%\end{strip}
% check that!
% \begin{align}
%   \label{fenchel_dual_ra2}
%   \cB(V)^\pfd(\pi_0, \gamma_0)
%   & = \begin{array}[t]{rl}
%         \inf\limits_{\va{\gamma}, \va{\lambda}, \va{\mu}}
%         & \left\langle \va\lambda, \va{d} \right\rangle
%         + \nesp{ V^\pfd(\va\mu - A^\top \va\lambda, \va\gamma) }
%         + \sup_{x_0} \ \pi_0^\top x_0 
%         - \left\langle B^\top \va\lambda, x_0\right\rangle \\
%         \textrm{s.t.}
%         & \va\gamma \va{c}^\top + T^\top \va\lambda \geq 0 \\
%         & \va\gamma \in \gamma_0 \cP, \va\mu \geq 0.
%   \end{array} \\
%   & = \begin{array}[t]{rl}
%         \inf\limits_{\va{\gamma}, \va{\lambda}, \va{\mu}}
%         & \left\langle \va\lambda, \va{d} \right\rangle
%         + \nesp{ V^\pfd(\va\mu - A^\top \va\lambda, \va\gamma) } \\
%         \textrm{s.t.}
%         & \nesp{B^\top \va\lambda} = \pi_0 \\
%         & \va\gamma \va{c}^\top + T^\top \va\lambda \geq 0 \\
%         & \va\gamma \in \gamma_0\cP, \va\mu \geq 0.
%   \end{array} \\
%   & = \begin{array}[t]{rl}
%         \inf\limits_{\va{\gamma}, \va{\lambda}, \va{\pi}}
%         & \left\langle \va\lambda, \va{d} \right\rangle
%         + \nesp{ V^\pfd(\va\pi, \va\gamma) } \\
%         \textrm{s.t.}
%         & \nesp{B^\top \va\lambda} = \pi_0 \\
%         & \va\pi + A^\top \va\lambda \geq 0 \\
%         & \va\gamma \va{c}^\top + T^\top \va\lambda \geq 0 \\
%         & \va\gamma \in \gamma_0\cP.
%   \end{array}
% \end{align}
This equation defines a risk-\emph{neutral} LBO $\cB^\pfd$
that takes a \emph{homogeneous} recourse function $V^\pfd$
and returns another homogeneous convex function of the same dimension.
We call this operator the \emph{projective dual Bellman operator} associated to $\cB$.
%This motivates the following
% \begin{definition}
%   The projective dual Bellman operator associated to $\cB$ is the LBO
%   \begin{equation}
%     \label{projective_dual_lbo}
%     \cB^\pfd(V^\pfd)(\pi_0, \gamma_0)
%     =  \begin{array}[t]{rll} 
%   \inf\limits_{ {\gamma}, {\lambda}, {\pi}}&
%     \sum_{j\in[J]} \lambda_j^\top d_j  + V^\pfd(\pi_j, \gamma_j)
%         \\[.3ex]
%     \textrm{s.t.} \quad
%      & \gamma \in \gamma_0 \cQ  \\
%     & \gamma_j c_j + T_j^\top \lambda_j \geq 0 & \forall j  \\
%     &  \sum_j B_j^\top\lambda_j = \pi_0 \\
%     & \pi_j + A_j^\top \lambda_j \geq 0  & \forall j .
%   \end{array}
%   \end{equation}
% \end{definition}

\new{Comparing~\eqref{eq:dualRCR} and~\eqref{fenchel_dual_ra2},
we notice the decomposition is not done at the same time-step for all variables:
in the first one, $\zeta$ is a single variable,
relaxing the incoming dual state constraint;
whereas in the second, it relaxes the outgoing dual state constraint.
Substituting $\pi_j = \psi_j + \zeta_j - \mu_j$, we obtain}
the following proposition,
linking the
coperspectives of the primal value functions
with the value functions of the dual problem.
\begin{proposition}
For $t \in[T]$, if the dual value function $D_t$ is defined by~\eqref{eq:dualRCR},
and $V_t$ is the primal value function defined by~\eqref{eq:DP_ra_primal}
then
\[
  D_t(\pi_t, \gamma_t)
  = - \inf_{\substack{\zeta_t + \psi_t \geq \pi_t \\ \zeta_t \geq 0}}
        \bar{x_t}^\top\zeta_t + V_t^\pfd(\psi_t, \gamma_t).
\]
In particular, $D_t$ is a concave, positively homogeneous,
one-sided Lipschitz regularization of $V_t^\pfd$.

Further, the value of primal Problem~\eqref{pb:ra-mlsp} is
% \begin{equation}
%   \label{eq:zeroth}
 $ \sup_{\pi_0}  \quad \pi_0 ^\top x_{0} + D_0(\pi_0, 1)$.
% \end{equation}
\end{proposition}

This proposition paves the way to a dual SDDP algorithm.
Indeed, it was shown in~\cite{leclere2020exact}
that SDDP can be applied to any sequence of functions
linked through linear Bellman operators (LBO) like $\cB^\pfd$.

% vim:set spelllang=en:

\section{Examples}
\label{sec:examples}

In this section, we provide an algorithm, in the lineage of SDDP,
for the risk-averse dual problem given by the recursion~\eqref{projective_dual_lbo_compacted}.
Then, we close with one numerical example from a real-world problem.
A more comprehensive discussion on the algorithm,
including implementation details,
can be found in the appendix.
There, one will also find further results
on the application of our algorithm.

\subsection{A dual risk-averse algorithm}
\label{sec:algorithm}
The recursion of (perspective) value functions $D_t$
given by~\eqref{projective_dual_lbo_compacted}
can be solved by recursively constructing piecewise linear (upper) approximations,
which we call $\mathfrak{D}_t$.
As usual, one needs to ensure that
the domain of the state variables $\pi_t$ and $\gamma_t$ remains bounded.
Since all $\gamma_t$ remain in $[0,1]$,
we only need bounds for $\pi_t$, which we assume are given by the user
as the
Lipschitz constants $L_t$
for the primal value functions $V_t$.
% "as discussed in section 2.4"?
In our experiments, the Lipschitz constant estimation was not critical:
Increasing $L_t$ by a factor $10$ or $100$ had a negligible impact after $50$ iterations,
as can be seen in section~C of the companion.
Moreover, one needs a starting upper bound for $\Dfrak_t$.
These can be obtained, for example,
choosing $\pi_t = 0$ and $\gamma_t = 1$,
and constructing cuts from $t = T-1$ back to $t = t_0$.

%As we observed in section~\ref{sub:RCR},
The first stage problem, corresponding to $t = t_0$, is slightly different.
It is obtained as the fusion of the ``zero-th stage''
containing $\pi_{n_0}$ as a decision variable,
and the first stage in~\eqref{projective_dual_lbo_compacted}.
Furthermore, since $x_{n_0}$ is fixed,
there's no corresponding slack variables $\mu_{n_0}$ and $\zeta_{n_0}$,
so it must satisfy
\begin{equation}
  \label{pi_0}
  \sum_{j \in [J_0]} B_j^\top \lambda_j = \pi_0 .
\end{equation}

With this, we can now present how one can perform
Bellman iterations on the recursion defined by~\eqref{projective_dual_lbo_compacted}
to obtain convergence.
We highlight the following differences with the primal SDDP:
\begin{itemize}
  \item Computing \new{$D_t(\pi, \gamma)$}
    cannot be decomposed by realization of $\va \omega_{t}$
    due to the coupling constraint
    $\zeta + \sum_{j \in [J_t]} B_j^\top \lambda_j \geq \pi_t$.
    In particular, the forward pass is as demanding as the backward pass,
    and yields cuts.
    Furthermore, we have one next-state variable
    per possible realization of $\va \omega_t$, which means that,
    when adding a single cut to the approximation of $D_{t+1}$,
    we are adding $J_t$ constraints.
  \item In the forward step, we choose the realization $j$ according to a
    (smoothed) ``importance sampling'' procedure, with weight
    $\gamma_j + \varepsilon$.
  \item By homogeneity, we normalize the state variables $(\pi_j, \gamma_j)$
    that will be used in the next step of the forward pass
    to have $\gamma_{t+1} = 1$, unless we are in a branch where $\gamma_t = 0$.
    This has had a positive impact in the numerical stability of the algorithm.
  \item Finally, by remark~\ref{rk:beta0}, we ensure that, for every cut,
    its parameter $\beta$ is always zero.
\end{itemize}

\begin{algorithm}
  \DontPrintSemicolon
  \KwData{upper bounds~$\Dfrak_t^0 \geq D_t$ and bounds~$L_t$ for $|\pi_t|$}
  % \KwData{maximum number of iterations $N$ and smoothing parameter $\varepsilon > 0$}
  \KwResult{upper bound on the value of \eqref{projective_dual_lbo_compacted}}
  \BlankLine
  %Set initial iteration $k=1$\;
  \For(){$k = 0$ \KwTo $N$}{
    Solve the first stage problem to obtain $\pi_{0}$, and set $\gamma_{0} = 1$ \;
    \lIf{$k == N$}
    {
      Return upper bound \;
    }

    \For(\tcp*[f]{forward pass}){$t = 0$ \KwTo $T-1$}{
      Solve problem~\eqref{projective_dual_lbo_compacted} with $\Dfrak_{t+1}^k$ instead of $D_{t+1}$\;
      Compute a cut for $D_t$ using the optimal multipliers for $\pi_t$ and $\gamma_t$ \;
      Choose a branch $\hat{\jmath}$ according to probabilities $\gamma_j + \varepsilon$ \;
      \eIf{$\gamma_{\hat{\jmath}} > 0$}{
        Set $\pi_{t+1} \leftarrow \pi_{\hat{\jmath}}/\gamma_{\hat{\jmath}}$, and $\gamma_{t+1} \leftarrow 1$
      }{
        Set $\pi_{t+1} \leftarrow \pi_{\hat{\jmath}}$, and $\gamma_{t+1} \leftarrow 0$
      }
    }
  }
  %Return the upper bound\;
\caption{Dual Risk-Averse SDDP}\label{alg:D-RA-SDDP}
\end{algorithm}

Naturally, one can couple this algorithm with (say) SDDP running on the primal.
This keeps track of both upper and lower bounds,
therefore allowing to stop based on a prescribed \emph{tolerance},
instead of just a maximum number of iterations as described above.

Let us close this section with two remarks.
First, even if this algorithm uses only forward passes,
one could use backward passes for computing cuts,
as in the classical SDDP algorithm.
This would require solving approximately twice the number of optimization problems,
but would include in the backward pass the updated value function,
which could potentially speed up the convergence of the algorithm.
Furthermore, this algorithm is easily amenable to standard cut-selection techniques,
which can be useful to reduce the computational burden of each iteration.

% vim:set spelllang=en:

\subsection{Numerical experiments}
\label{sec:numerical}
%!TEX root=../main.tex

\new{We present here a numerical example.
Further details and other results are given in the companion,
and the implementation in julia, along with other examples, can be found at
\url{https://github.com/bfpc/DualSDDP.jl}.}

This example comes from the Brazilian Hydrothermal Energy planning problem,
where the reservoirs and hydro dams are aggregated into $4$ subsystems,
and there is a 5th node in the network, as an interconnection.
Therefore, it contains $4$ state variables (the stored energy in each reservoir),
$9$ equality constraints for the dynamics
($4$ for the states, and $5$ for demand in each node),
and a total of $164$ control variables,
accounting for hydro and thermal energy produced,
and energy exchange among the nodes in the system.
The uncertainty at each time step is the inflow for each aggregated reservoir,
and is different for each time step, corresponding to different months of the year.

For this example, we take $12$ stages and $82$ inflow realizations per stage
(thus $82^{12}$ scenarios).
We have natural bounds for every state variable,
given by the reservoirs' limits,
and control variables (power output, line capacities, \ldots).
The risk measure considered was a combination of expectation and AV@R,
given by $\new{\beta}\mathbb{E} + \new{(1-\beta)}\text{AV@R}_{\new{\alpha}}$.
\new{In this problem, the highest marginal cost is given by load shedding,
which yields estimates for the Lipschitz constants we use.}

\begin{figure}[H]
  \centering
  \includegraphics[width=0.375\textwidth]{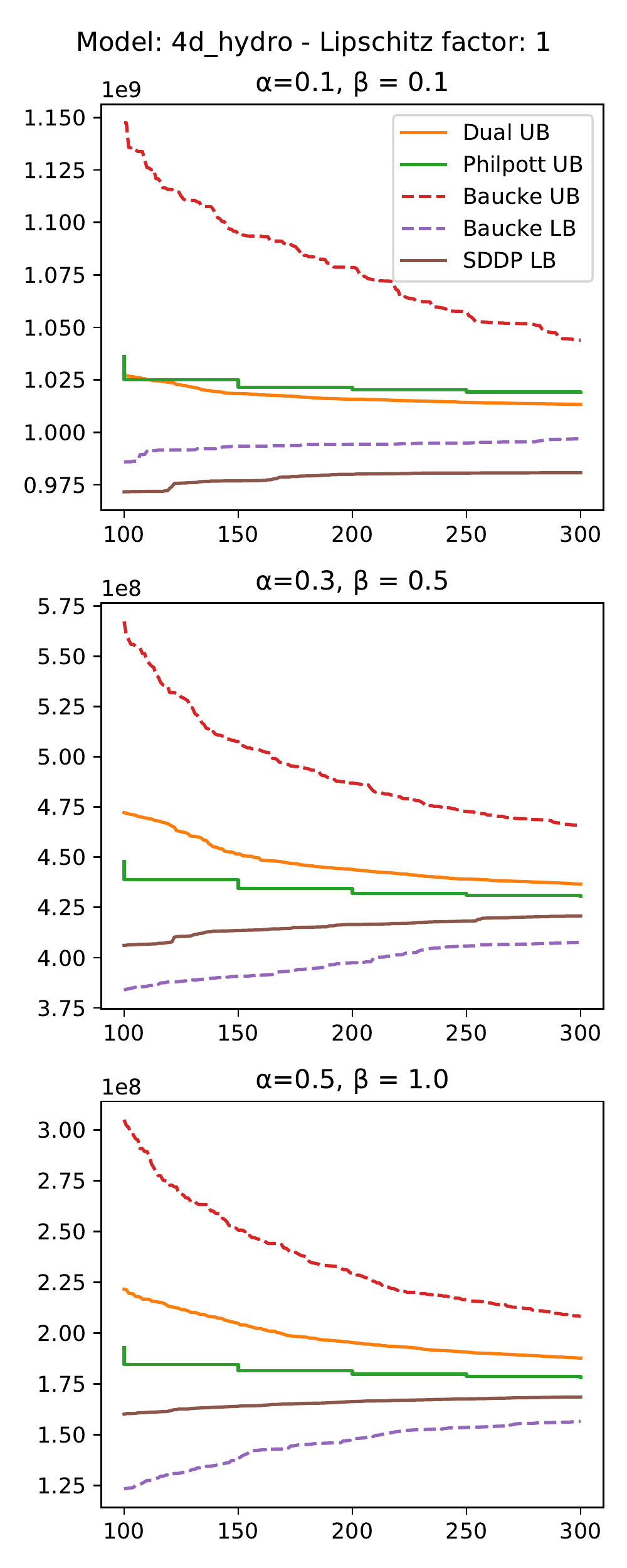}
  \caption{Bounds evolution for hydrothermal problem.}
  \label{fig:4d_diag}
\end{figure}

\new{In Figure~\ref{fig:4d_diag}, we present the evolution of the bounds obtained by the primal SDDP, our dual SDDP algorithm, as well the one shot backward bounds of \cite{philpott2013solving} (Philpott UB),
 computed every $50$ iterations based on the trajectories from primal SDDP, and the upper and lower bounds provided by the problem-child method of \cite{baucke2017deterministic} (Baucke UB / LB). 
 This is done for various level of risk aversion. 
 Note that, on this problem, the dual upper bound always outperform the problem-child method.
 It also slightly beat the primal one-shot upper bound in the most risk-averse case.
 This is also observed on the other numerical experiments available at \url{https://github.com/bfpc/DualSDDP.jl}.
 }

\new{Finally, we noticed that each iteration of the dual is between 30 and 15 times slower than primal iteration, being larger for higher branching sizes.
%Further, this relative difference decreases as iterations progress.
}

\begin{table}[H]
  \centering
  {\small 
  \begin{tabular}[h]{c | c c c}
    \# branches
    & P-SDDP & D-SDDP & Problem Child \\
    \hline 
    10 & 0.023 & 0.166 & 0.109\\
20 & 0.054 & 0.523 & 0.224\\
40 & 0.113 & 2.366 & 0.402\\
80 & 0.274 & 5.739 & 0.813
  \end{tabular}
  \caption{Single iteration time in sec (around $it=100$)}
  }
\end{table}

This is expected, since each problem in the dual formulation
includes all inflow realizations and a linking constraint among all of them,
whereas the primal problem also allows decomposing each time step
in separate problems for each branch.

% vim:set spelllang=en:

\paragraph{Acknowledgements}
We would like to thank the editor and an anonymous referee
for their suggestions which improved the presentation of our results.
We also thank Lucas Merabet for his comments.

The first author is partly supported by project COPPETEC-23145.
The second author benefited from the support of
FMJH-PGMO and from EDF.

\bibliographystyle{alpha}
\bibliography{biblio}

\appendix
\section{Dual SDDP algorithm}
Algorithm~\ref{alg:dual_SDDP_detailed} presents the details of the dual dynamic programming algorithm
used to solve the regularized dual problem in the recursion~(13),
which we recall here for ease of reference:

\begin{subequations}
  \label{dual_rec}
  \begin{alignat}{5}
    \label{eq:D_T}
    D_T(\pi_T, \gamma_T) & \omit\rlap{${} = - \bar{x}_T^\top \max(\pi_T, 0)$} \\
    \label{projective_dual_lbo_appendix}
    D_t(\pi_t, \gamma_t)
    & = & \sup\limits_{\zeta, \gamma_j, \lambda_j, \pi_j, \xi_j}
        & - \bar{x}_t^\top \zeta +
            \sum\limits_{j\in[J]} - d_j^\top \lambda_j
            - \bar{y}_{t+1}^\top \xi_j + D_{t+1}(\pi_j, \gamma_j) \hspace{-5em} \\[.3ex]
    && \textrm{s.t.} \quad
         & \gamma \in \gamma_t \cQ \label{cst:coupling_gamma} \\
    &&   & \zeta + \sum_j B_j^\top\lambda_j \geq \pi_t \label{cst:coupling_pi}\\
    &&   & \pi_j + A_j^\top \lambda_j = 0  & \forall j \in [J_t] \\
    &&   & \gamma_j c_j + T_j^\top \lambda_j + \xi_j \geq 0 & \forall j \in [J_t] \\
    &&   & |\pi_j|\leq \gamma_j L_{t+1} & \forall j \in [J_t] \\
    &&   & \zeta \geq 0, \xi_j \geq 0. \label{cst:last_rec}
  \end{alignat}
\end{subequations}
where the recursion defined by equations
\eqref{projective_dual_lbo_appendix}--\eqref{cst:last_rec}
applies for all $t \in \{0, 1, \ldots, T-1\}$.

Recall that $\bar x_t$ and $\bar y_{t+1}$ are upper bounds
on the norm of the state and control variables,
and $L_{t+1}$ is a Lipschitz constant on the primal value function $V_t$.
Note that $\zeta$ and $\xi_j$ can be interpreted as slack variables in the dual with exact penalization given by the bounds on the primal variables.

%%% Former discussion with linearization
% The regularized backwards recursion starts from time $T$ with
% \begin{equation}
%   \label{eq:D_T}
%   D_T(\pi_T, \gamma_T) = - \bar{x}_T^\top \max(\pi_T, 0) \,,
% \end{equation}
% which can be exactly modeled in the linear problem
% for the last recursion, linking $D_{T-1}$ to $D_T$,
% using extra variables.
% More precisely, we substitute $D_T(\pi_j, \gamma_j)$
% in the objective of~\eqref{projective_dual_lbo_appendix}
% by $- \bar{x}_T^\top \zeta_j$,
% and add the pair of constraints
% $\zeta_j \geq 0$, $\zeta_j \geq \pi_j$
% for every branch $j \in [J_T]$.

\begin{algorithm}[H]
  \DontPrintSemicolon

  \KwData{Valid upper bounds $\D_t^0$ for the dual functions $D_t$}
  \KwData{Lipschitz constants $L_t$ bounding the dual states $\pi_t$}
  \KwPar{Number of iterations $N$}
  \KwPar{Tolerance $tol$ for small probabilities}
  \KwPar{Smoothing constant $\varepsilon$ for scenario sampling}
  \KwResult{Upper bound for problem~(1)}
  \KwResult{Improved upper bounds $\D_t^N$ for the dual functions $D_t$}
  \BlankLine

  Set $\text{UpperBound} \leftarrow +\infty$ \;
  \For{$k = 0$ \KwTo $N$}{
    \tcp{Compute initial state and update upper bound}
    Solve $\sup_{\pi_0} \pi_0^\top x_0 + \D_0^k(\pi_0, 1)$,
    save the optimal value $d$, optimal state $\pi_0$
    and set $\gamma_0 \leftarrow 1$ \nllabel{inistate} \;

    Set $\text{UpperBound} \leftarrow d$ \;
    \lIf{$k == N$}
    {
      \Return UpperBound, $\D_t^N$
    }

    \BlankLine
    \tcp{Forward Pass}
    \For{$t = 0$ \KwTo $T-1$}{
      Solve the optimization problem: \nllabel{Qtilde}
      \begin{equation}
        \tag{*}
        \label{eq:Qtilde}
        D_t(\pi_t, \gamma_t)
        =  \begin{array}[t]{rll}
          \sup\limits_{\zeta, \gamma_j, \lambda_j, \pi_j, \xi_j} &
          - \bar{x}_t^\top \zeta +
          \sum\limits_{j \in [J_t]} \left( - d_j^\top \lambda_j
          - \bar{y}_{t+1}^\top \xi_j + \D^k_{t+1}(\pi_j, \gamma_j) \right) \\[.3ex]
          \textrm{s.t.} \quad
          & \gamma \in \gamma_t \cQ  \\
          & \zeta + \sum_j B_j^\top\lambda_j \geq \pi_t \\
          & \pi_j + A_j^\top \lambda_j = 0  & \forall j \in [J_t] \\
          & \gamma_j c_j + T_j^\top \lambda_j + \xi_j \geq 0 & \forall j \in [J_t] \\
          & |\pi_j|\leq \gamma_j L_{t+1} & \forall j \in [J_t] \\
          & \zeta \geq 0, \xi_j \geq 0
        \end{array}
      \end{equation}
      save the optimal variables $\pi_j$, $\gamma_j$,
      and the optimal value $v_t$ \;

      \BlankLine
      \tcp{Construct a valid cut and add it to $\D_t$}
      Set $x_t \leftarrow \text{the dual multiplier of the constraint 
      $\zeta + \sum_j B_j^\top \lambda_j \geq \pi_t$}$\;
      Set $z_t \leftarrow \text{the dual multiplier for
      $\gamma \in \gamma_t \mathcal{Q}$}$\;

      \tcp{Sanity check}
      \lIf{$v_t \not\approx x_t^\top \pi_t + z_t \cdot \gamma_t$}
      {
        Warn numerical instability
      }
      \BlankLine

      Let $C(\pi, \gamma) := x_t^\top \pi + z_t \cdot \gamma$\;

      Set $\D_t^{k+1} \leftarrow \min(\D_t^k, C)$ \;

      \BlankLine
      \tcp{Prepare for next stage}
      Choose the next branch $\jh$ among all $J_t$ branches
      with probability proportional to $\gamma_j + \epsilon$
      \nllabel{next_state} \;

      \uIf(\tcp*[f]{Normalize state}){$\jh > \text{tol}$}
      {
        Set $\pi_{t+1} \leftarrow \pi_{\jh} / \gamma_{\jh}$,
        and $\gamma_{t+1} \leftarrow 1$\nllabel{line:normalization} \;
      } \Else(\tcp*[f]{Round down}) {
        Set $\pi_{t+1} \leftarrow \pi_{\jh}$,
        and $\gamma_{t+1} \leftarrow 0$ \nllabel{line:rounding} \;
      }
    }

  }
  \caption{Dual SDDP algorithm, for risk-averse problems}
  \label{alg:dual_SDDP_detailed}
\end{algorithm}

\newpage 

The dual algorithm is initialized with upper-approximations $\D_t^0$
of the dual value functions $D_t$,
which must be guaranteed upper bounds.
One possibility is to compute them from the costs $c_t$
and the bounds $\bar{y}_t$ of the control variables of the primal problem,
since the stage costs are at most
$\max_{0 \leq \va{y}_t \leq \bar{y}_t} \va{c}_t^\top \va{y}_t$.
% do we give more details?
Another possibility consists in doing a backward pass,
as described in Remark~\ref{rk:backward} below,
on any admissible dual trajectory (e.g. $(\pi_t=0, \gamma_t=1)$ for all $t$).

Line~\ref{inistate} of Algorithm~\ref{alg:dual_SDDP_detailed} computes an initial dual state $\pi_0$
given the current approximation $\D_0^k$. 
Indeed, from proposition~5, we know that the primal problem has optimal value
\begin{equation}
  \label{eq:marg_cost}
  \sup_{\pi_0} \pi_0^\top x_0 + D_0(\pi_0, 1) \,,
\end{equation}
where $x_0$ is the primal initial state.
 This dual state can (and often does) change between iterations.

% upon convergence, it will yield the marginal cost of the initial state.

From this initial state, Problem~\eqref{eq:Qtilde} in line~\ref{Qtilde}
is analogous to equation~\eqref{projective_dual_lbo_appendix}, with $\D_{t+1}^k$ in place of $D_{t+1}$.
In an LP implementation, $\D_{t+1}^k(\pi_j, \gamma_j)$ can be represented
through a hypographical variable $z_j$ for each scenario $j$,
and cuts $C$ for $\D_{t+1}$ become linear constraints of the form
\begin{equation}
  z_j \leq C^\kappa(\pi_j, \gamma_j), \qquad \kappa \in[k].
\end{equation}
Solving Problem~\eqref{eq:Qtilde} yields solutions $(\pi_j,\gamma_j)_{j\in[J]}$
corresponding to the outgoing states for all realizations of $\va \xi_{t+1}$.

Line~\ref{next_state} randomly selects the next state,
in a way that each branch as a positive probability to be chosen,
with a preference towards the scenarios which most contribute to the value function.
More precisely, the probability of choosing a branch $j$
is proportional to $\gamma_j + \epsilon$,
where $\epsilon$ is a small positive number
and $\va\gamma$ the current change-of-measure.
Note that for some risk measures, like the AV@R
(but not strict combinations of AV@R and Expectation),
the current change of measure could attribute $0$ probability to some realizations,
preventing exploration, and thus convergence of the dynamic programming algorithm.

By homogeneity of the value functions,
the probabilities $\gamma_t$ are normalized in line~\ref{line:normalization} at each stage,
and $\pi_t$ is normalized accordingly.
We have observed that this usually improves the numerical stability of the algorithm.
Indeed, the value of $\gamma_t$
is the (current) risk-adjusted probability of the stage-$t$ scenario,
which decreases as $t$ increases.
Since most solvers have both a relative and an absolute tolerance,
the homogeneity of the stage problems with respect to $(\pi_t, \gamma_t)$
might result in a very large relative error of the algorithm
when $\gamma_t$ becomes too small.
Finally, in line~\ref{line:rounding}, the probabilities $\gamma_t$ are rounded down to $0$ if they are too small.

%TODO
After performing $N$ iterations, the algorithm stops,
returning the current best upper bound for Problem~1,
and the current piecewise linear approximations $\D_t^N$
of the dual functions $D_t$. 

\begin{remark}[Backward pass]
  \label{rk:backward}
In the primal SDDP algorithm (risk-averse or not), the stage problem can be decomposed in $J_t$ subproblems, one for each realization of $\va \xi_{t+1}$. 
In particular, the optimal next-state for realization $j\in[J_t]$ is given by solving a problem independent of other possible realizations of $\va \xi_{t+1}$.
However, computing a cut requires solving a problem that depends on all realizations of $\va \xi_{t+1}$.
Thus, standard implementations of primal SDDP have a forward phase, to determine trajectories,
and a backward phase, to compute cuts;
the latter is slower, solving $[J_t]$ more problems at stage $t$.

In the dual formulation, this decomposition is no longer possible
due to the coupling constraints \eqref{cst:coupling_gamma} and~\eqref{cst:coupling_pi}.
In particular, to determine the optimal next-state value for a given realization $j \in [J_t]$,
one needs to solve a problem that depends on all realizations of $\va \xi_{t+1}$.
Thus, computing a dual trajectory also provides all the information needed to compute a cut. 
This is why Algorithm~\ref{alg:dual_SDDP_detailed} only has a forward phase. 

Naturally, it is also possible to add cuts in a backward fashion,
    which would need then to solve a problem similar to equation~\eqref{eq:Qtilde},
    but with an extra cut, using $\D_{t+1}^{k+1}$.
    This speeds up the information flow back to the first stage,
    at the cost of (approximately) doubling the time per iteration.
    This might be especially useful in the first few iterations to replace the initial, user-given, upper bound.
\end{remark}

\section{Detailed description of the numerical experiments}

The numerical example we used comes from the Brazilian Hydrothermal Energy planning problem.
In its long-term formulation,
the reservoirs and hydro dams are aggregated into $4$ subsystems,
Southeast, South, Northeast and North.
Each subsystem also corresponds to a region with an associated total energy demand.
Long-distance transmission lines connect the South with the Southeast,
Southeast with Northeast,
and an extra interconnection node (modeled as a $5^\text{th}$ subsystem),
to the North, Northeast and Southeast subsystems.
In each subsystem, the demand for energy in each month, $d_{s,t}$,
is supposed to be known;
the demand of subsystem 5 is zero.
Not satisfying this demand with thermal or hydro-generation
and exchanges with another subsystem,
leads to energy curtailment,
as described in~\eqref{cst:demand_satisfaction}.

For simplicity, this model considers energy equivalents for water volumes,
so the stored volumes are represented by $x_{s,t}$,
the equivalent energy in the reservoir of subsystem $s$
at the end of stage $t$ (and the beginning of stage $t+1$).
For system $s$, the hydro generation during stage $t$ is given by $h_{s,t}$,
the (equivalent energy) inflow by $\text{inflow}_{s,t}$,
and (equivalent energy) spillage by $\text{spill}_{s,t}$,
resulting in the dynamic equation~\eqref{cst:dyn}.
Constraints~\eqref{cst:dam_capacity} to~\eqref{cst:exchange_bound} represent
physical bounds on hydro storage, hydro production, thermal production and exchanges.
The spillage is akin to a slack variable, and therefore positive
as enforced by~\eqref{cst:spill}.
Remaining constraints define four ranges for energy curtailment.

Thermal power plants are represented individually,
each with its own minimum and maximum generation limits,
$\underline{g}_j$, $\overline{g}_j$,
as well as costs per MWh $c_j$.
Each thermal plant is located in a given subsystem $s$,
and the set $T_s$ collects the indices $j$ of thermal plants
in subsystem $s$.
If demand is not met, curtailment has increasing costs
$CD_k$ for $k = \{1,2,3,4\}$,
corresponding to curtailment below 5\%, 10\%, 20\% or 100\%
of the demand $d_{s,t}$ of the subsystem.

Therefore, the (primal) dynamic programming recursion becomes:
\begin{subequations}
  \label{hydroproblem}
  \begin{align}
  V_t(x_{t-1}) =
  \min        \quad & \rho\Big[ \sum_j c_j g_{j,t} + \sum_k \sum_s CD_k cur_{s,k,t}
                          + \text{spill\_pen} \sum_s \text{spill}_{s,t} \nonumber \\
              &    \qquad + \sum_{s,s'} \text{xch\_pen}_{s,s'} ex_{s,s',t}
                          + V_{t+1}(x_t) \Big] \\[0.5ex]
  \text{s.t.} \quad
              & d_{s,t} = h_{s,t} + \sum_{j \in T_s} g_{j,t} + \sum_k cur_{s,k,t} +
              \sum_{s'} ex_{s',s,t} - ex_{s,s',t},  & \forall t, \forall s \label{cst:demand_satisfaction}\\
              & x_{s,t} = x_{s,t-1} + \text{inflow}_{s,t} - h_{s,t} - \text{spill}_{s,t} & \forall t, \forall s \label{cst:dyn}\\
              & 0 \leq x_{s,t} \leq \overline{x}_s& \forall t, \forall s \label{cst:dam_capacity}\\
              & 0 \leq h_{s,t} \leq \overline{h}_s &\forall t, \forall s \label{cst:hydro_capacity}\\
              & \underline{g}_j \leq g_{j,t} \leq \overline{g}_j& \forall t, \forall j \label{cst:thermal_bound}\\
              & 0 \leq ex_{s,s',t} \leq \overline{ex}_{s,s'}& \forall t, \forall s, \forall s'  \label{cst:exchange_bound}\\
              & 0 \leq \text{spill}_{s,t}& \forall t, \forall s \label{cst:spill} \\
              & 0 \leq cur_{s,1,t} \leq 5\% \cdot d_{s,t}  & \forall t, \forall s \label{cst:curtailment_1}\\
              & 0 \leq cur_{s,2,t} \leq 5\% \cdot d_{s,t}  & \forall t, \forall s \label{cst:curtailment_2}\\
              & 0 \leq cur_{s,3,t} \leq 10\% \cdot d_{s,t} & \forall t, \forall s \label{cst:curtailment_3}\\
              & 0 \leq cur_{s,4,t} \leq 80\% \cdot d_{s,t} & \forall t, \forall s \label{cst:curtailment_4}.
  \end{align}
\end{subequations}

% where the constraints are, in order:
% \begin{enumerate}
%   \item The hydro balance equation, in energy-equivalent, on the reservoir $s$;
%   \item The demand satisfaction for subsystem $s$, from hydro and thermal generation,
%     curtailment in 4 levels and exchanges between subsystems;
%   \item Bounds for stored energy in each equivalent reservoir;
%   \item Bounds for hydro and thermal generation;
%   \item Bounds for each curtailment level;
%   \item Bounds for energy exchange from subsystem $s$ to subsystem $s'$,
%     which are zero if there's no connection, and are not necessarily symmetrical;
%   \item Bounds for (energy equivalent) spillage in the hydro plants.
% \end{enumerate}
The stage costs include thermal generation costs,
and curtailment costs for every level and subsystem.
Moreover, it includes penalties for both energy spillage and exchange.

The problem instances we solve consider uncertainties on the inflows only.
We take the historical inflows for each month as scenarios,
which are then sampled independently along the planning horizon.
This amounts to 82 realizations per stage,
corresponding to the years 1931--2012, inclusive.

Data, such as variable bounds and unit costs for the example
we deal with can be found at the supplementary file \texttt{data.jl}.
A further supplementary file \texttt{demand.jl}
contains the series of demands, for each subsystem,
along the stages.
The historical series of inflows we use
can be found in the last supplementary file, \texttt{eafs.npz}.

\begin{sloppypar}
A complete setup, parsing the data and building the corresponding
matrices for the dual recursion can be found
at~\url{https://github.com/bfpc/DualSDDP.jl/blob/91a50a9c9eb16db6acc4a046e4471c9737cd01a1/examples/4d_hydro/}.
\end{sloppypar}

\section{Impact of Lipschitz estimate on convergence}

We performed two experiments to assess the impact of
providing a larger Lipschitz constant than the true one.
In order to do so, we used algorithm~\ref{alg:dual_SDDP_detailed}
with a tight Lipschitz constant, then a $10$ times larger one,
and finally a $100$ times larger one.
We assessed different combinations of risk-aversion,
and compared the evolution of the upper bounds
to the best lower bound found with the primal SDDP.

\begin{figure}[ht]
  \centering
  \includegraphics[width=0.9\textwidth]{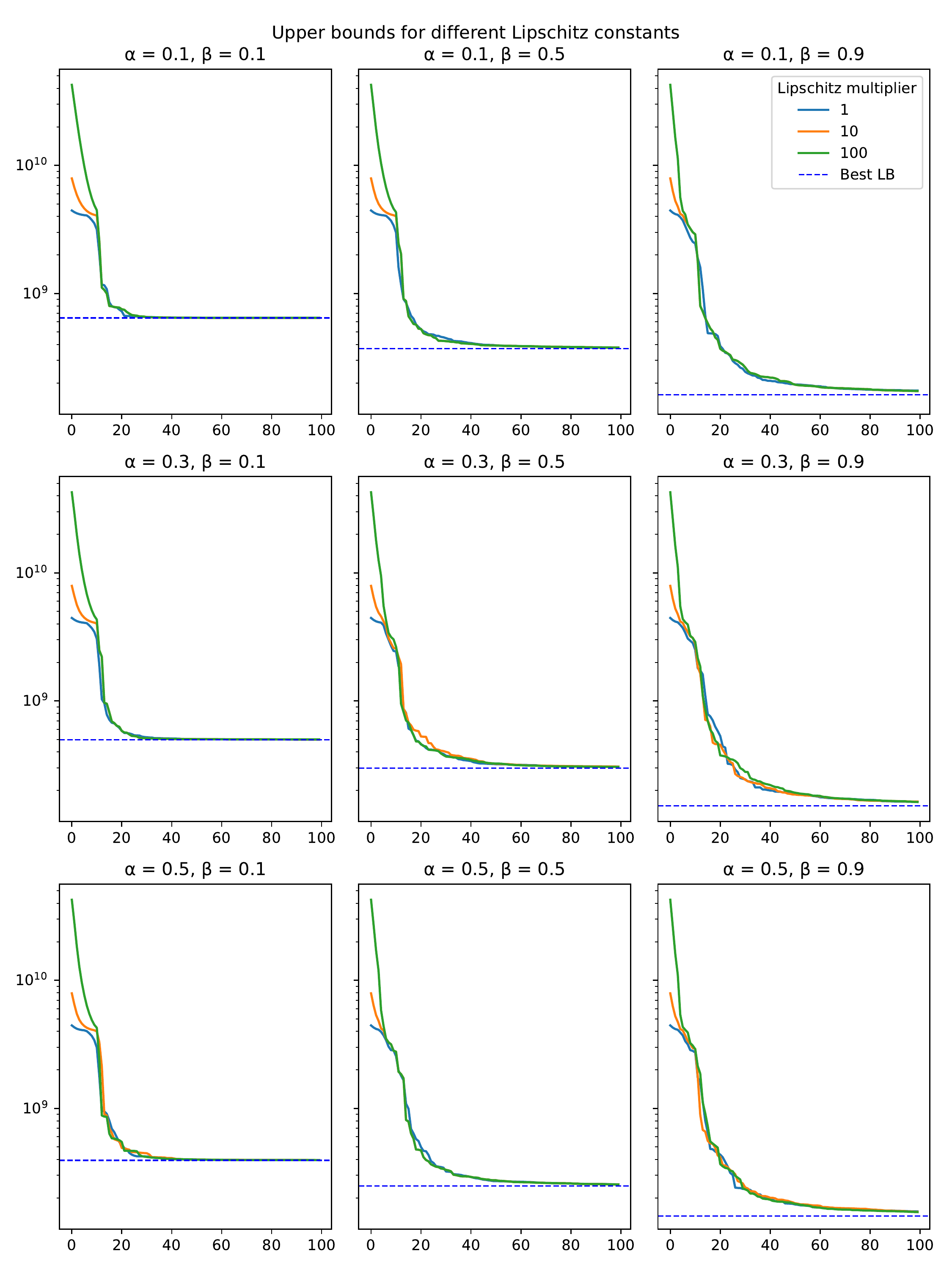}
  \caption{Upper bounds on the small 2-reservoir problem.}
  \label{fig:Lip-2d}
\end{figure}

The first graph, in Figure~\ref{fig:Lip-2d},
corresponds to a simplified hydrothermal problem,
given by the same recursion~\eqref{hydroproblem},
but with only three thermal units, two reservoirs,
and one interconnection between the two corresponding subsystems.
We notice that the initial estimates are larger
for larger Lipschitz estimates,
but after some iterations the impact of a
worse Lipschitz estimate is negligible.

\begin{figure}[ht]
  \centering
  \includegraphics[width=0.9\textwidth]{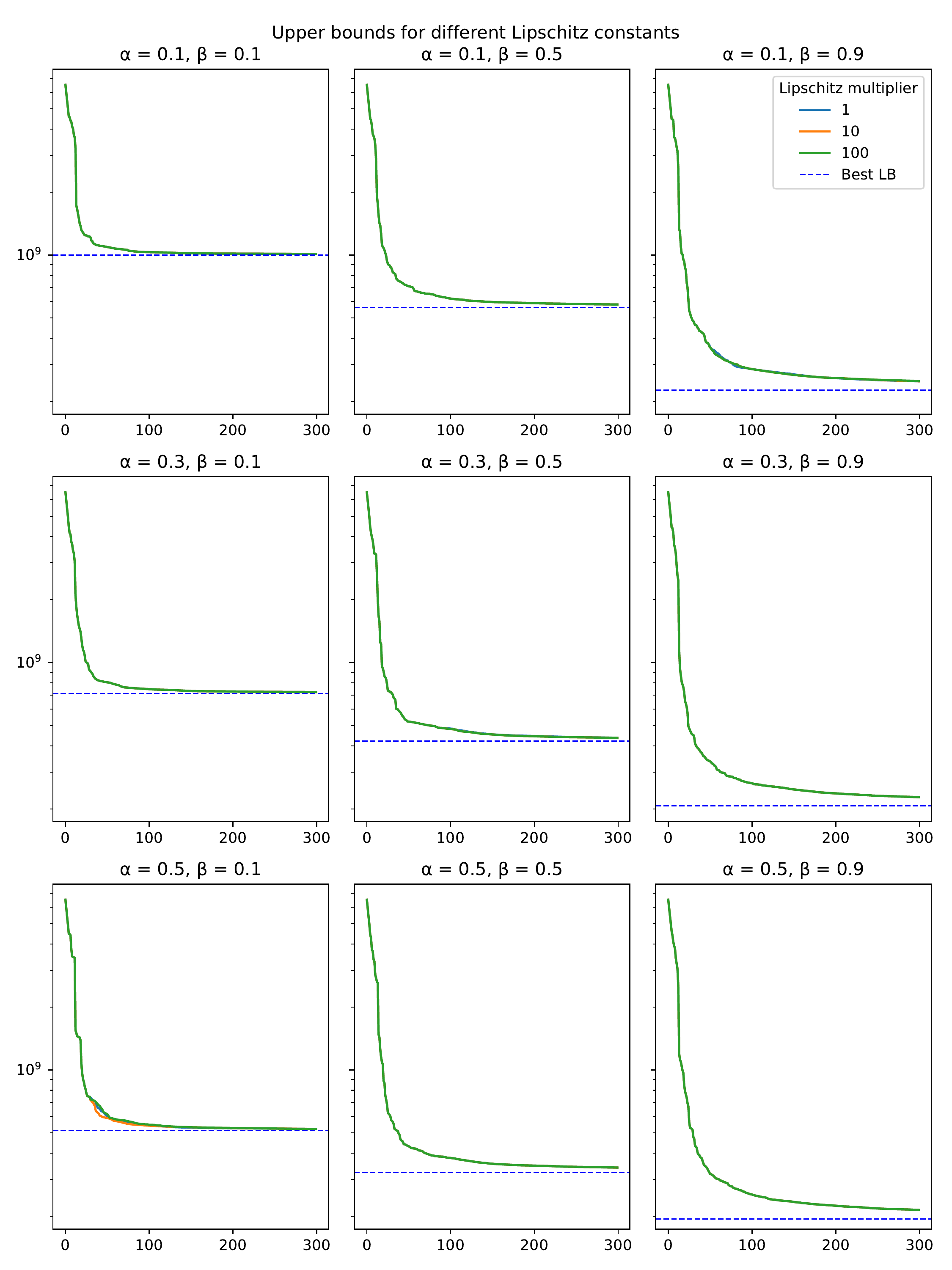}
  \caption{Upper bounds on the 4-reservoir problem.}
  \label{fig:Lip-4d}
\end{figure}

The second one, in Figure~\ref{fig:Lip-4d},
corresponds to the larger 4-reservoir setting
of the previous section.
There, we remark a much lower sensitivity of the bounds
with respect to the candidate Lipschitz constant.
For example, the gaps at the 100th iteration
in the case of $\alpha = 0.3$ and $\beta = 0.5$
are, respectively, 14.95, 14.72 and~14.72
for factors 1, 10 and~100,
which is such a small difference that it is not visible in the figure.

For completeness, we report the relative gaps, in \%,
for both experiments in Table~\ref{table_gaps_lip},
for several intermediate iterations.

\begin{table}[ht]
  \centering
\begin{tabular}{|cc|rrrrr|rrrr|}
  \hline
                  &         & \multicolumn{5}{c|}{2 reservoir}& \multicolumn{4}{c|}{4 reservoir} \\
                  &         & \multicolumn{5}{c|}{Iteration}  & \multicolumn{4}{c|}{Iteration}  \\
 $(\alpha,\beta)$ &  Factor &
                1 &      10 &       20 &       50 &      100 &      10  &      100 &      200 &      300 \\ \hline
 \multirow{ 3}{*}{(0.10, 0.10)}
 &    1 &   590.4 &   448.6 &    16.32 &     0.42 &     0.19 &   304.04 &     3.66 &     1.94 &     1.61 \\
 &   10 &  1129.1 &   538.4 &    20.58 &     0.37 &     0.21 &   304.04 &     3.66 &     1.94 &     1.61 \\
 &  100 &  6515.6 &   659.4 &    20.58 &     0.37 &     0.21 &   304.04 &     3.66 &     1.94 &     1.61 \\ \hline
 \multirow{ 3}{*}{(0.10, 0.50)}
 &    1 &  1097.4 &   817.5 &    47.49 &     6.58 &     1.58 &   542.00 &    10.23 &     4.86 &     3.30 \\
 &   10 &  2031.6 &  1000.5 &    42.94 &     5.73 &     1.84 &   542.00 &    10.23 &     4.86 &     3.30 \\
 &  100 & 11373.8 &  1138.4 &    42.94 &     5.73 &     1.84 &   542.00 &    10.23 &     4.86 &     3.30 \\ \hline
 \multirow{ 3}{*}{(0.10, 0.90)}
 &    1 &  2643.2 &  1462.4 &   187.70 &    20.90 &     7.53 &  1436.20 &    26.41 &    14.24 &    10.31 \\
 &   10 &  4783.5 &  1757.4 &   168.36 &    21.39 &     6.73 &  1436.20 &    26.45 &    14.37 &    10.33 \\
 &  100 & 26186.4 &  1757.4 &   168.36 &    21.39 &     6.73 &  1436.20 &    26.45 &    14.42 &    10.66 \\ \hline
 \multirow{ 3}{*}{(0.30, 0.10)}
 &    1 &   794.4 &   596.9 &    23.97 &     1.25 &     0.65 &   383.48 &     5.16 &     2.42 &     1.88 \\
 &   10 &  1492.3 &   722.2 &    26.51 &     1.16 &     0.69 &   383.48 &     5.16 &     2.42 &     1.88 \\
 &  100 &  8470.8 &   828.2 &    26.51 &     1.16 &     0.69 &   383.48 &     5.16 &     2.42 &     1.88 \\ \hline
 \multirow{ 3}{*}{(0.30, 0.50)}
 &    1 &  1385.2 &   721.1 &    61.54 &     7.44 &     2.21 &   687.45 &    14.95 &     5.90 &     3.85 \\
 &   10 &  2544.0 &   765.4 &    94.13 &     8.12 &     2.54 &   687.45 &    14.72 &     5.64 &     3.83 \\
 &  100 & 14131.6 &   911.3 &    60.85 &     8.74 &     2.11 &   687.45 &    14.72 &     5.64 &     3.83 \\ \hline
 \multirow{ 3}{*}{(0.30, 0.90)}
 &    1 &  2833.4 &  1771.2 &   283.57 &    24.23 &     7.30 &  1505.73 &    28.37 &    14.88 &    10.06 \\
 &   10 &  5122.1 &  1954.3 &   202.89 &    22.78 &     8.02 &  1505.73 &    28.37 &    14.88 &    10.06 \\
 &  100 & 28009.1 &  1943.0 &   210.18 &    27.93 &     7.25 &  1505.73 &    28.37 &    14.88 &    10.06 \\ \hline
 \multirow{ 3}{*}{(0.50, 0.10)}
 &    1 &  1032.7 &   768.6 &    39.17 &     1.46 &     0.57 &   574.96 &     6.06 &     2.46 &     1.72 \\
 &   10 &  1916.5 &   938.9 &    41.30 &     1.55 &     0.57 &   574.96 &     5.50 &     2.67 &     1.75 \\
 &  100 & 10754.4 &  1046.1 &    43.92 &     1.65 &     0.59 &   574.96 &     6.67 &     2.74 &     1.72 \\ \hline
 \multirow{ 3}{*}{(0.50, 0.50)}
 &    1 &  1695.5 &  1047.3 &   124.57 &     9.26 &     2.91 &   920.73 &    17.26 &     7.67 &     5.38 \\
 &   10 &  3096.3 &  1040.5 &    92.07 &    11.03 &     2.92 &   920.73 &    17.26 &     7.67 &     5.38 \\
 &  100 & 17104.7 &  1040.5 &    92.07 &    11.03 &     2.92 &   920.73 &    17.26 &     7.67 &     5.45 \\ \hline
 \multirow{ 3}{*}{(0.50, 0.90)}
 &    1 &  2976.4 &  1846.5 &   210.09 &    24.36 &     8.41 &  1659.86 &    31.14 &    15.66 &    10.35 \\
 &   10 &  5376.6 &  1945.9 &   196.51 &    27.87 &     7.61 &  1659.86 &    31.14 &    15.66 &    10.35 \\
 &  100 & 29379.1 &  2039.5 &   241.80 &    25.54 &     7.31 &  1659.86 &    31.14 &    15.66 &    10.35 \\ \hline
\end{tabular}
  \caption{Relative gaps (\%) for the 2-reservoir and 4-reservoir problems,
  for different factors corresponding to overestimating the Lipschitz constant.}
  \label{table_gaps_lip}
\end{table}

% vim:set spelllang=en:

\end{document}